\documentclass[12pt]{amsart}

\usepackage{amsmath,xspace,amssymb,mathrsfs}
\usepackage{color}

\input xy
\xyoption{all}
\xyoption{2cell}
\UseAllTwocells
\CompileMatrices

\newcommand{\Spec}{\operatorname{Spec}}
\renewcommand{\phi}{\varphi}

\newcommand{\rn}{\operatorname{rank}}
\newcommand{\hm}{\operatorname{Hom}}

\newcommand{\Ker}{\operatorname{Ker}}

\newcommand{\Ima}{\operatorname{Im}}

\newcommand{\MA}{\operatorname{Max}}

\newcommand{\Min}{\operatorname{Min}}
\newcommand{\Ann}{\operatorname{Ann}}

\newcommand{\Supp}{\operatorname{Supp}}

\newtheorem{proposition}{Proposition}[section]
\newtheorem{lemma}[proposition]{Lemma}

\newtheorem{corollary}[proposition]{Corollary}
\newtheorem{theorem}[proposition]{Theorem}

\theoremstyle{definition}
\newtheorem{definition}[proposition]{Definition}
\newtheorem{example}[proposition]{Example}

\newtheorem{remark}[proposition]{Remark}

\begin{document}

\title[Projectivity of finitely generated flat modules]{On the projectivity of finitely generated flat modules}

\author[A. Tarizadeh]{Abolfazl Tarizadeh}
\address{Department of Mathematics, Faculty of Basic Sciences, University of Maragheh \\
P. O. Box 55136-553, Maragheh, Iran.
 }
\email{ebulfez1978@gmail.com}

\date{}

\begin{abstract} In this paper, the projectivity of a finitely generated flat module of a commutative ring is studied through its exterior powers and invariant factors and then various new results are obtained. Specially, the related results of Endo, Vasconcelos, Wiegand, Cox-Rush and Puninski-Rothmaler on the projectivity of finitely generated flat modules are generalized.

\end{abstract}

\maketitle

\section{Introduction}

The main purpose of the present paper is to investigate the projectivity of finitely generated flat modules of a commutative ring. This topic has been the main subject of many articles in the literature over the years and it is still of current interest, see e.g. \cite{Cox-Pendleton}, \cite{Endo}, \cite{Facchini}, \cite{Jondrup}, \cite{Puninski-Rothmaler}, \cite{Vasconcelos}, \cite{Wiegand}. Note that in general there are finitely generated flat modules which are not projective, see Example \ref{example 1}, also see \cite[Tag 00NY]{Johan} as another example (note that our example is so simple than the cited one; it is also applicable for other purposes). The main motivation to investigate the projectivity of finitely generated flat modules essentially originates from the fact that ``every finitely generated flat module over a local ring is free'',  see Theorem \ref{lemma 8}. This result together with Theorem \ref{th 4} plays a major role in this paper. \\

In this paper, the projectivity of a finitely generated flat module of a commutative ring is studied through its exterior powers and invariant factors and then we obtain various new and interesting results. One of the features of this study is that some major results in the literature on the projectivity of finitely generated flat modules are generalized. In particular, Theorem \ref{th 6} generalizes \cite[Theorem 1]{Endo}, Theorem \ref{th 3} a little improves \cite[Theorem 2.1]{Vasconcelos}, Theorem \ref{coro 5} generalizes  \cite[Theorem 2]{Wiegand}, \cite[Proposition 2.3]{Cox-Pendleton}, it also generalizes \cite[Proposition 5.5 and Corollary 5.6]{Puninski-Rothmaler} in the commutative case and Corollary \ref{Corollary 19} generalizes \cite[Theorem 2]{Wiegand}. In fact, Theorem \ref{coro 5} can be viewed as a generalization of all of the above mentioned results. This theorem is one of the main novel contributions of this paper and has many non-trivial consequences. Theorems \ref{Corollary VII}, \ref{Corollary VIII} and \ref{Theorem I} are another interesting results of this paper. \\

For reading the present paper having a reasonable knowledge on the ``exterior powers of a module'' is necessary. In this paper, all rings are commutative. \\

\section{Preliminaries}

We need the following material in the next sections.  \\

\begin{lemma}\label{lemma 2} Let $M$ be a finitely generated $R-$module, let $I=\Ann_{R}(M)$ and let $S$ be a multiplicative subset of $R$. Then $S^{-1}I=\Ann_{S^{-1}R}(S^{-1}M)$.\\
\end{lemma}

{\bf Proof.} Easy. $\Box$ \\

A projective $R-$module is also called $R-$projective. We also use a similar terminology for free and flat modules. \\

Unlike the Kaplansky theorem \cite{Kaplansky} which states that every projective module over a local ring is free, but this is not true for flat modules. For example, the field of rationals $\mathbb{Q}$ is $\mathbb{Z}_{\mathfrak{p}}-$flat but it is not $\mathbb{Z}_{\mathfrak{p}}-$free where $\mathfrak{p}$ is a non-zero prime ideal of the ring of integers $\mathbb{Z}$. In spite of this, in the finite case we have the following interesting result which can be considered it as the analogue of the Kaplansky theorem for flat modules and plays a major role in this paper.\\

\begin{theorem}\label{lemma 8} Every finitely generated flat module over a local ring is free. \\
\end{theorem}

{\bf Proof.} See \cite[Theorem 7.10]{Matsumura} or \cite[Tag 00NZ]{Johan}. $\Box$ \\

\begin{lemma}\label{lemma 7} Let $\phi:R\rightarrow S$ be a morphism of rings and $M$ a finitely generated flat $R-$module. Then $\Ann_{S}(M\otimes_{R}S)=\Ann_{R}(M)S$. \\
\end{lemma}

{\bf Proof.} It is a local property and implies from Lemma \ref{lemma 2} and Theorem \ref{lemma 8}. $\Box$ \\

If $R\rightarrow S$ is a ring map and $M$ is an $R-$module then $\Lambda^{n}(M)\otimes_{R}S$ as $S-$module is canonically isomorphic to $\Lambda_{S}^{n}(M\otimes_{R}S)$. It is also well-known that if $M$ is a projective (resp. flat) $R-$module then for each natural number $n$,  $\Lambda^{n}(M)$ is a projective (resp. flat) $R-$module. Finally, if $M$ is a finitely generated $R-$module then $\Lambda^{n}(M)$ is a finitely generated $R-$module. \\

If $M$ is a $R-$module then the $n$-th invariant factor of $M$, denoted by $I_{n}(M)$, is defined as the annihilator of the $n$-th exterior power of  $M$, i.e., $I_{n}(M)=\Ann_{R}\big(\Lambda^{n}(M)\big)$. \\

\begin{remark}\label{Remark 040} If $M$ is a finitely generated flat $R-$module then Theorem \ref{lemma 8} leads us to a function $\psi:\Spec R\rightarrow\mathbb{N}=\{0,1,2,...\}$ which is defined as $\mathfrak{p}\rightsquigarrow\rn_{R_{\mathfrak{p}}}(M_{\mathfrak{p}})$. It is called the rank map of $M$. It is obvious that the rank map is continuous if and only if it is locally constant (i.e., for each prime ideal $\mathfrak{p}$ of $R$ then there exists an open neighborhood $U\subseteq\Spec(R)$ of that point such that $\rn_{R_{\mathfrak{q}}}(M_{\mathfrak{q}})
=\rn_{R_{\mathfrak{p}}}(M_{\mathfrak{p}})$ for all $\mathfrak{q}\in U$). It is well known that
$\Supp\big(\Lambda^{n}(M)\big)=\{\mathfrak{p}\in\Spec(R): \rn_{R_{\mathfrak{p}}}(M_{\mathfrak{p}})\geq n\}$. \\
\end{remark}

If $\phi:R\rightarrow S$ is a morphism of rings then the induced map $\Spec(S)\rightarrow\Spec(R)$ given by $\mathfrak{p}\rightsquigarrow\phi^{-1}(\mathfrak{p})$ is denoted by $\phi^{\ast}$ or by $\Spec(\phi)$. \\

The radical Jacobson of a ring $R$ is denoted by $\mathfrak{J}(R)$. \\

An ideal $I$ of a ring $R$ is called a pure ideal if the canonical ring map $R\rightarrow R/I$ is a flat ring map. Pure ideals are quite interesting and play an important role in commutative and non-commutative algebra (for instance, in classifying Gelfand rings and their dual rings). \\

\begin{theorem}\label{Remark 030} An ideal $I$ of a ring $R$ is a pure ideal if and only if $\Ann(f)+I=R$ for all $f\in I$. \\
\end{theorem}

{\bf Proof.} It is a local property and implies from Theorem \ref{lemma 8}. $\Box$ \\

\begin{corollary}\label{Corollary I} Let $M$ be a finitely generated flat $R-$module with the annihilator $I$. Then $I$ is a pure ideal.  \\
\end{corollary}

{\bf Proof.} It is a local property and implies from Theorem \ref{lemma 8} and Theorem \ref{Remark 030}. $\Box$ \\

\begin{lemma}\label{lemma 1} The annihilator of a finitely generated projective module is generated by an idempotent element.\\
\end{lemma}

{\bf Proof.} It is well known and implies from the fact that if $M$ is a projective $R-$module then $JM=M$ where $J=\big(\phi(m): \phi\in\hm_{R}(M,R), m\in M\big)$ is the trace ideal of $M$. $\Box$ \\

The following result is well known, see \cite[Chap. II, \S 5.2, Th\'{e}or\`{e}me 1]{Bourbaki}, \cite[Tag 00NX]{Johan} and \cite[Proposition 1.3]{Vasconcelos}. As a contribution to this result, we provide a new proof for the equivalency $\mathbf{(ii)}$. \\

\begin{theorem}\label{th 4} Let $M$ be a finitely generated flat $R-$module. Then the following are equivalent.\\
$\mathbf{(i)}$ $M$ is $R-$projective.\\
$\mathbf{(ii)}$ The invariant factors of $M$ are finitely generated ideals.\\
$\mathbf{(iii)}$ The rank map of $M$ is locally constant. \\
\end{theorem}

{\bf Proof.} $\textbf{(i)}\Rightarrow\textbf{(ii)}:$ It is well-known that $\Lambda^{n}(M)$ is a finitely generated projective $R-$module and so by Lemma \ref{lemma 1}, $I_{n}(M)$ is a principal ideal.\\
$\textbf{(ii)}\Rightarrow\textbf{(iii)}:$ It suffices to show that the rank map of $M$ is Zariski continuous. By Corollary \ref{Corollary I}, $I_{n}(M)$ is an idempotent ideal. Thus there exists some $a\in I_{n}(M)$ such that $(1-a)I_{n}(M)=0$. Clearly $a=a^{2}$ and $I_{n}(M)=Ra$. By Remark \ref{Remark 040}, $\psi^{-1}(\{n\})=\Supp N\cap\big(\Spec(R)\setminus\Supp N'\big)$ where $N=\Lambda^{n}(M)$ and $N'=\Lambda^{n+1}(M)$. But $\Supp N=D(1-a)$. Moreover, $\Supp N'=V\big(I_{n+1}(M)\big)$ since $N'$ is a finitely generated $R-$module. Therefore $\psi^{-1}(\{n\})$ is an open subset of $\Spec R$. \\
$\textbf{(iii)}\Rightarrow\textbf{(i)}:$ Apply Theorem \ref{lemma 8} and \cite[Tag 00NX]{Johan}. $\Box$ \\

In the following we give an example of a finitely generated flat module which is not projective. It should be noted that finding such examples of modules is not as easy as one may think at first.\\

\begin{example}\label{example 1} Let $R=\prod\limits_{i\geq1}A$ be an infinite product of copies of a non-zero ring $A$ and let $I=\bigoplus\limits_{i\geq1}A$ which is an ideal of $R$. If $f=(f_{i})\in I$ then there exists a finite subset $D$ of $\{1,2,3,...\}$ such that $f_{i}=0$ for all $i\in\{1,2,3,...\}\setminus D$. Clearly $f=fg$ where $g=(g_{i})\in I$ with $g_{i}=1$ for all $i\in D$ and $g_{i}=0$ for all $i\in\{1,2,3,...\}\setminus D$. Hence, $I$ is a pure ideal of $R$ (i.e., $R/I$ is a finitely generated flat $R-$module). If $R/I$ is $R-$projective then by Lemma \ref{lemma 1}, there exists a sequence $e=(e_{i})\in R$ such that $I=Re$. Thus there exists a finite subset $E$ of $\{1,2,3,...\}$ such that $e_{i}=0$ for all $i\in\{1,2,3,...\}\setminus E$. Clearly $\{1,2,3,...\}\setminus E\neq\emptyset$.
Pick some $k\in\{1,2,3,...\}\setminus E$. There is some $r=(r_{i})\in R$ such that $(\delta_{i,k})_{i\geq1}=re$ where $\delta_{i,k}$ is the Kronecker delta. In particular, $1_{A}=r_{k}e_{k}=r_{k}0_{A}=0_{A}$. This is a contradiction. Therefore $R/I$ is not $R-$projective. \\
\end{example}

\section{Projectivity- main results}

Throughout this section, $M$ is a finitely generated flat $R-$module. \\

The following technical result generalizes \cite[Theorem 1]{Endo}. \\

\begin{theorem}\label{th 6} Let $R\rightarrow S$ be an injective ring map.
Then $M$ is $R-$projective if and only if $M\otimes_{R}S$ is $S-$projective.\\
\end{theorem}

{\bf Proof.} Let $M\otimes_{R}S$ be $S-$projective. Without loss of generality, we may assume that $R\subseteq S$ is an extension of rings. First we shall prove that $I=\Ann_{R}(M)$ is a principal ideal.
By Lemma \ref{lemma 7}, $IS=L$ where $L=\Ann_{S}(N)$ and $N=M\otimes_{R}S$. Hence by Lemma \ref{lemma 1}, there is an idempotent $e\in S$ such that $IS=Se$. Let $J=S(1-e)\cap R$. Clearly $IJ=0$. We have $I+J=R$. If not, then there exists a prime ideal $\mathfrak{p}$ of $R$ such that $I+J\subseteq\mathfrak{p}$. Thus, by Theorem \ref{lemma 8}, $I_{\mathfrak{p}}=0$. Therefore the extension of $IS$ under the canonical map $S\rightarrow S\otimes_{R}R_{\mathfrak{p}}$ is zero. Thus there exists an element $s\in R\setminus\mathfrak{p}$ such that $se=0$ and so $s=s(1-e)$. Hence $s\in J$. But this is a contradiction. Therefore $I+J=R$. It follows that there is an element $c\in I$ such that $c=c^{2}$ and $I=Rc$. Now, let $n\geq1$. We have $\Lambda^{n}(M)$ is a finitely generated flat $R-$module. Moreover, $\Lambda^{n}(M)\otimes_{R}S$ is $S-$projective because it is canonically isomorphic to $\Lambda^{n}_{S}(M\otimes_{R}S)$. Thus, by what we have proved above, $I_{n}(M)$ is a principal ideal. Hence, by Theorem \ref{th 4}, $M$ is $R-$projective. The reverse is easy and well known. $\Box$ \\

The following result is also technical and generalizes \cite[Theorem 2.1]{Vasconcelos}.\\

\begin{theorem}\label{th 3} Let $J$ be an ideal of $R$ which is contained in the radical Jacobson of $R$. If $M/JM$ is $R/J-$projective then $M$ is $R-$projective.\\
\end{theorem}

{\bf Proof.} First we shall prove that $I=\Ann_{R}(M)$ is a principal ideal.
By Lemma \ref{lemma 7}, $L=I+J$ where $L=\Ann_{R}(M/JM)$. Also, by Lemma \ref{lemma 1},  $\Ann_{R/J}(M/JM)=L/J$ is a principal ideal. This implies that $I=Rx+I\cap J$ for some $x\in R$ since $L/J=I+J/J$ is canonically isomorphic to $I/(I\cap J)$. But $I=Rx$. Because let $\mathfrak{m}$ be a maximal ideal of $R$. By Theorem \ref{lemma 8}, $I_{\mathfrak{m}}$ is either the whole localization or the zero ideal. If $I_{\mathfrak{m}}=0$ then $(Rx)_{\mathfrak{m}}=0$ since $Rx\subseteq I$. But if $I_{\mathfrak{m}}=R_{\mathfrak{m}}$ then $I$ is not contained in $\mathfrak{m}$. Thus $Rx$ is also not contained in $\mathfrak{m}$  since $I\cap J\subseteq J\subseteq\mathfrak{m}$. Hence $(Rx)_{\mathfrak{m}}=R_{\mathfrak{m}}$. Therefore $I=Rx$.
Now let $n\geq1$ and let $N=\Lambda^{n}(M)$. Then $N/JN$ is $R/J-$projective. Because, $N/JN$ as $R/J-$module is canonically isomorphic to $\Lambda^{n}_{R/J}(M/JM)$ and $\Lambda^{n}_{R/J}(M/JM)$ is $R/J-$projective. But $N$ is a finitely generated flat $R-$module. Therefore, by what we have proved above, $I_{n}(M)=\Ann_{R}(N)$ is a principal ideal. Thus the invariant factors of $M$ are finitely generated ideals and so by Theorem \ref{th 4}, $M$ is $R-$projective.
$\Box$ \\

Motivated by the Grothendieck's relative point of view, then we obtain the following result which (beside Theorems \ref{lemma 8} and \ref{th 4}) is one of the most powerful results on the projectivity of finitely generated flat modules. \\

\begin{theorem}\label{coro 5} Let $\phi:R\rightarrow S$ be a ring map whose kernel is contained in the radical Jacobson of $R$. Then $M$ is $R-$projective iff $M\otimes_{R}S$ is $S-$projective. \\
\end{theorem}

{\bf Proof.} Let $M\otimes_{R}S$ be $S-$projective. Clearly $M/JM$ is a finitely generated flat $R/J-$module and $M/JM\otimes_{R/J}S\simeq M\otimes_{R}S$ is $S-$projective where $J=\Ker\phi$. Moreover $R/J$ can be viewed as a subring of $S$ via $\phi$.
Therefore, by Theorem \ref{th 6}, $M/JM$ is $R/J-$projective. Then by applying Theorem \ref{th 3}, we get that $M$ is $R-$projective. The reverse is easy and well known. $\Box$ \\

The above theorem has many consequences. \\

Recall that a ring $R$ is called an S-ring (``S'' referes to Sakhajev) if every finitely generated flat $R-$module is $R-$projective. \\

\begin{corollary}\label{Corollary XII} Let $\phi:R\rightarrow S$ be a ring map whose kernel is contained in the radical Jacobson of $R$. If $S$ is an S-ring then $R$ is as well.\\
\end{corollary}

{\bf Proof.} If $M$ is a finitely generated flat $R-$module then $M\otimes_{R}S$ is a finitely generated flat $S-$module and so, by the hypothesis, it is $S-$projective. Therefore by Theorem \ref{coro 5}, $M$ is $R-$projective. $\Box$ \\

\begin{corollary} If there exists a prime ideal $\mathfrak{p}$ of a ring $R$ such that the kernel of the canonical ring map $\pi:R\rightarrow R_{\mathfrak{p}}$ is contained in the radical Jacobson of $R$, then $R$ is an S-ring. \\
\end{corollary}

{\bf Proof.} By Theorem \ref{lemma 8}, every local ring is an S-ring, then apply
Corollary \ref{Corollary XII}. $\Box$ \\

\begin{corollary} If the radical Jacobson of a ring $R$ contains a prime ideal of $R$, then $R$ is an S-ring. $\Box$ \\
\end{corollary}

\begin{remark} Let $S$ be a subset of a ring $R$. The polynomial ring $R[x_{s} : s\in S]$ modulo $I$ is denoted by $S^{(-1)}R$ where the ideal $I$ is generated by elements of the form $sx_{s}^{2}-x_{s}$ and $s^{2}x_{s}-s$ with $s\in S$. We call $S^{(-1)}R$ the pointwise localization of $R$ with respect to $S$. Amongst them, the pointwise localization of $R$ with respect to itself, namely $R^{(-1)}R$, has more interesting properties; for further information please consult with \cite{Olivier}. Note that Wiegand \cite{Wiegand} utilizes the notation $\widehat{R}$ instead of $R^{(-1)}R$. Clearly $\eta(s)=\eta(s)^{2}(x_{s}+I)$ and $x_{s}+I=\eta(s)(x_{s}+I)^{2}$ where $\eta:R\rightarrow S^{(-1)}R$ is the canonical map and the pair $(S^{(-1)}R, \eta)$ satisfies in the following universal property: ``for each such pair $(A,\phi)$, i.e. $\phi:R\rightarrow A$ is a ring map and for each $s\in S$ there is some $c\in A$ such that $\phi(s)=\phi(s)^{2}c$ and $c=\phi(s)c^{2}$, then there exists a unique ring map $\psi:S^{(-1)}R\rightarrow A$ such that $\phi=\psi\circ\eta$". Now let $\mathfrak{p}$ be a prime ideal of $R$ and consider the canonical map $\pi:R\rightarrow\kappa(\mathfrak{p})$ where $\kappa(\mathfrak{p})$ is the residue field of $R$ at $\mathfrak{p}$. By the above universal property, there is a (unique) ring map $\psi:S^{(-1)}R\rightarrow\kappa(\mathfrak{p})$ such that $\pi=\psi\circ\eta$. Thus $\eta$ induces a surjection between the corresponding spectra. This, in particular, implies that the kernel of $\eta$ is contained in the nil-radical of $R$. Using this, then the following result generalizes \cite[Theorem 2]{Wiegand}.\\
\end{remark}

\begin{corollary}\label{Corollary 19} If there exists a subset $S$ of $R$ such that $M\otimes_{R}S^{(-1)}R$ is $S^{(-1)}R-$projective then $M$ is $R-$projective. \\
\end{corollary}

{\bf Proof.} It is an immediate consequence of Theorem \ref{coro 5}. $\Box$ \\

In what follows we get some new results essentially based on the referee's excellent comments. \\

\begin{corollary} Let $\phi:R\rightarrow S$ be a morphism of rings such that the induced map $\phi^{\ast}$ has the dense image. Then $M$ is $R-$projective iff $M\otimes_{R}S$ is $S-$projective. \\
\end{corollary}

{\bf Proof.} It is an immediate consequence of Theorem \ref{coro 5}, because from $\overline{\Ima\phi^{\ast}}=\Spec(R)$ we get that the kernel of $\phi$ is contained in $\mathfrak{J}(R)$, the radical Jacobson of $R$. $\Box$ \\

\begin{lemma}\label{Lemma III} Let $(R_{k},\mathfrak{m}_{k})$ be a family of local rings. Then the kernel of the canonical ring map $R:=\prod\limits_{k}R_{k}\rightarrow\prod\limits_{k}R_{k}/\mathfrak{m}_{k}$ is contained in the radical Jacobson of $R$. \\
\end{lemma}

{\bf Proof.} If the sequence $x=(x_{k})\in R$ is a member of the kernel then $x_{k}\in\mathfrak{m}_{k}$ for all $k$. To prove the assertion it suffices to show $1+xy$ is invertible in $R$ for all $y=(y_{k})\in R$. For each $k$, there exists some $z_{k}\in R_{k}$ such that $(1+x_{k}y_{k})z_{k}=1$ because $R_{k}$ is a local ring. It follows that
$(1+xy)z=1$ where $z=(z_{k})$. $\Box$ \\

\begin{corollary}\label{Corollary IX} Let $X\subseteq\Spec(R)$ be a subset. Then the following are equivalent. \\
$\mathbf{(i)}$ $M\otimes_{R}S$ is $S-$projective where $S=\prod\limits_{\mathfrak{p}\in X}\kappa(\mathfrak{p})$. \\
$\mathbf{(ii)}$ $M\otimes_{R}S'$ is $S'-$projective where $S'=\prod\limits_{\mathfrak{p}\in X}R/\mathfrak{p}$. \\
$\mathbf{(iii)}$ $M\otimes_{R}S''$ is $S''-$projective where $S''=\prod\limits_{\mathfrak{p}\in X}R_{\mathfrak{p}}$. \\
If moreover, $\bigcap\limits_{\mathfrak{p}\in X}\mathfrak{p}\subseteq\mathfrak{J}(R)$ then the above statements are equivalent with the following. \\
$\mathbf{(iv)}$ $M$ is $R-$projective. \\
\end{corollary}

{\bf Proof.} It is an immediate consequence of Theorem \ref{coro 5}. $\Box$ \\

The subsets $\Min(R)$ and $\MA(R)$ are typical examples which satisfy in the hypothesis of Corollary \ref{Corollary IX}. \\

\begin{corollary}\label{Corollary VI} Consider the following commutative diagram of rings:  $$\xymatrix{R\ar[r]\ar[d]&S'\ar[d]\\S\ar[r]^{\phi}&T}$$ in which the kernel of $\phi$ is contained in the radical Jacobson of $S$. If $M\otimes_{R}S'$ is $S'-$projective then $M\otimes_{R}S$ is $S-$projective.  \\
\end{corollary}

{\bf Proof.} If $M\otimes_{R}S'$ is $S'-$projective then it is easy to see that: $$(M\otimes_{R}S')\otimes_{S'}T\simeq M\otimes_{R}T\simeq(M\otimes_{R}S)\otimes_{S}T$$ is $T-$projective. But $M\otimes_{R}S$ is a finitely generated flat $S-$module. Therefore by Theorem \ref{coro 5}, $M\otimes_{R}S$ is $S-$projective. $\Box$ \\

\begin{definition} If $X$ is a subset of $\Spec(R)$ then we call $\bigcup\limits_{\mathfrak{p}\in X}V(\mathfrak{p})$ the \emph{specialization cone} of $X$ and it is denoted by $X_{s}$ or by $\mathcal{Z}(X)$. Dually, we call $\bigcup\limits_{\mathfrak{p}\in X}\Lambda(\mathfrak{p})$ the \emph{generalization cone} of $X$ and it is denoted by $X_{g}$ or by $\mathcal{F}(X)$ where $\Lambda(\mathfrak{p})=\{\mathfrak{q}\in\Spec(R):\mathfrak{q}\subseteq
\mathfrak{p}\}$. \\
\end{definition}

\begin{theorem}\label{Corollary VII} Let $X\subseteq\Spec(R)$ be a subset. Put $S:=\prod\limits_{\mathfrak{p}\in X}\kappa(\mathfrak{p})$ and $S':=\prod\limits_{\mathfrak{p}\in X_{s}}\kappa(\mathfrak{p})$. Then $M\otimes_{R}S$ is $S-$projective iff $M\otimes_{R}S'$ is $S'-$projective. \\
\end{theorem}

{\bf Proof.} Consider the canonical injective ring map $T=\prod\limits_{\mathfrak{p}
\in X}R/\mathfrak{p}\rightarrow S$. Then by Theorem \ref{th 6}, $M\otimes_{R}S\simeq(M\otimes_{R}T)\otimes_{T}S$ is $S-$projective iff $M\otimes_{R}T$ is $T-$projective. By the axiom of choice, we obtain a function $\sigma:X_{s}\rightarrow X$ such that $\sigma(\mathfrak{p})\subseteq\mathfrak{p}$ for all $\mathfrak{p}\in X_{s}$ and $\sigma(\mathfrak{p})=\mathfrak{p}$ for all $\mathfrak{p}\in X$. For each $\mathfrak{p}\in X$, consider the canonical injective ring map $R/\mathfrak{p}\rightarrow\prod\limits_{\mathfrak{q}
\in\sigma^{-1}(\mathfrak{p})}R/\mathfrak{q}$. Then we get the canonical injective ring map $T\rightarrow T'=\prod\limits_{\mathfrak{p}
\in X_{s}}R/\mathfrak{p}$. Again by Theorem \ref{th 6}, $M\otimes_{R}T$ is $T-$projective iff $M\otimes_{R}T'$ is $T'-$projective. Similarly above, by applying Theorem \ref{th 6} to the canonical injective ring map $T'\rightarrow S'$, we get that $M\otimes_{R}T'$ is $T'-$projective iff $M\otimes_{R}S'$ is $S'-$projective. $\Box$ \\

\begin{corollary}\label{Corollary X} Let $\mathfrak{p}$ be a prime ideal of a ring $R$ and put $S:=\prod\limits_{\mathfrak{q}\in V(\mathfrak{p})}\kappa(\mathfrak{q})$. Then $M\otimes_{R}S$ is $S-$projective. \\
\end{corollary}

{\bf Proof.} It is an immediate consequence of Theorem \ref{Corollary VII}. $\Box$ \\

\begin{theorem}\label{Corollary VIII} Let $X\subseteq\Spec(R)$ be a subset. Put $S:=\prod\limits_{\mathfrak{p}\in X}\kappa(\mathfrak{p})$ and $S':=\prod\limits_{\mathfrak{p}\in X_{g}}\kappa(\mathfrak{p})$. Then $M\otimes_{R}S$ is $S-$projective iff $M\otimes_{R}S'$ is $S'-$projective. \\
\end{theorem}

{\bf Proof.} The kernel of the canonical ring map $T=\prod\limits_{\mathfrak{p}\in X}R_{\mathfrak{p}}\rightarrow S$ is contained in the radical Jacobson of $T$, see Lemma \ref{Lemma III}. Therefore by Theorem \ref{coro 5}, $M\otimes_{R}S\simeq(M\otimes_{R}T)\otimes_{T}S$ is $S-$projective if and only if $M\otimes_{R}T$ is $T-$projective. By the axiom of choice, there exists a function $\sigma:X_{g}\rightarrow X$ such that $\mathfrak{p}\subseteq\sigma(\mathfrak{p})$ for all $\mathfrak{p}\in X_{s}$ and $\sigma(\mathfrak{p})=\mathfrak{p}$ for all $\mathfrak{p}\in X$. For each $\mathfrak{p}\in X$, consider the canonical injective ring map $R_{\mathfrak{p}}\rightarrow\prod\limits_{\mathfrak{q}
\in\sigma^{-1}(\mathfrak{p})}R_{\mathfrak{q}}$. Then we get the canonical injective ring map $T\rightarrow T'=\prod\limits_{\mathfrak{p}
\in X_{s}}R_{\mathfrak{p}}$. Thus by Theorem \ref{th 6}, $M\otimes_{R}T$ is $T-$projective iff $M\otimes_{R}T'$ is $T'-$projective. Again by Lemma \ref{Lemma III}, the kernel of the canonical ring map $T'\rightarrow S'$ is contained in the radical Jacobson of $T'$. Hence by Theorem \ref{coro 5}, $M\otimes_{R}T'$ is $T'-$projective if and only if $M\otimes_{R}S'$ is $S'-$projective. $\Box$ \\

\begin{corollary}\label{Corollary XI} Let $\mathfrak{p}$ be a prime ideal of a ring $R$ and put $S:=\prod\limits_{\mathfrak{q}\in\Lambda(\mathfrak{p})}
\kappa(\mathfrak{q})$. Then $M\otimes_{R}S$ is $S-$projective. \\
\end{corollary}

{\bf Proof.} It is an immediate consequence of Theorem \ref{Corollary VIII}. $\Box$ \\

\begin{corollary} Let $X\subseteq Y\subseteq\Spec(R)$ be two subsets such that either $Y\subseteq X_{s}$ or $Y\subseteq X_{g}$. Put $S:=\prod\limits_{\mathfrak{p}\in X}\kappa(\mathfrak{p})$ and $S':=\prod\limits_{\mathfrak{p}\in Y}\kappa(\mathfrak{p})$. Then $M\otimes_{R}S$ is $S-$projective iff $M\otimes_{R}S'$ is $S'-$projective. \\
\end{corollary}

{\bf Proof.} If $Y\subseteq X_{s}$ then $X_{s}=Y_{s}$, so apply Theorem \ref{Corollary VII} in this case. But if $Y\subseteq X_{g}$ then $X_{g}=Y_{g}$, and so apply Theorem \ref{Corollary VIII}. $\Box$ \\

\begin{lemma}\label{Lemma IV} Let $\phi:R\rightarrow S$ be a morphism of rings and put $S':=\prod\limits_{\mathfrak{p}\in\Ima\phi^{\ast}}\kappa(\mathfrak{p})$. Then $M\otimes_{R}S$ is $S-$projective if and only if $M\otimes_{R}S'$ is $S'-$projective. \\
\end{lemma}

{\bf Proof.} If $\mathfrak{p}\in\Ima\phi^{\ast}$ and
$\mathfrak{q}\in(\phi^{\ast})^{-1}(\mathfrak{p})$ then we have the canonical ring map $\kappa(\mathfrak{p})\rightarrow\kappa(\mathfrak{q})$ which is injective since every ring map from a field into a non-zero ring is injective. Then we get the canonical injective ring map $\kappa(\mathfrak{p})\rightarrow\prod\limits_
{\mathfrak{q}\in(\phi^{\ast})^{-1}(\mathfrak{p})}\kappa(\mathfrak{q})$. So we get the canonical injective ring map $S'\rightarrow T=\prod\limits_
{\mathfrak{q}\in\Spec(S)}\kappa(\mathfrak{q})$ which fits in the following commutative diagram: $$\xymatrix{R\ar[r]^{\phi}\ar[d]&S\ar[d]\\S'\ar[r]&T}$$
where the unnamed arrows are the canonical morphisms. It is easy to see that the kernel of the canonical morphism $S\rightarrow T$ is the nil-radical of $S$ which is contained in the radical Jacobson of $S$. Then the assertion is deduced by twice using of Corollary \ref{Corollary VI}. $\Box$ \\

Given a subset $X\subseteq\Spec(R)$, denote $X^{(1)}:=(X_{s})_{g}$ and $X_{(1)}:=(X_{g})_{s}$, and inductively $X^{(n)}:=(X^{(n-1)})^{(1)}$ and $X_{(n)}:=(X_{(n-1)})_{(1)}$. Note that in general, $X^{(n)}\neq X_{(n)}$. For example, if $X=\{2\mathbb{Z}\}\subseteq\Spec(\mathbb{Z})$ then $X^{(1)}=\{0,2\mathbb{Z}\}$ but $X_{(1)}=\Spec(\mathbb{Z})$. \\

\begin{theorem}\label{Theorem I} Let $\phi:R\rightarrow S$ be a morphism of rings and $X=\Ima\phi^{\ast}$. Assume there exists some $n\geq1$ such that $\bigcap\limits_{\mathfrak{p}\in X^{(n)}}\mathfrak{p}$ is contained in the radical Jacobson of $R$. If $M\otimes_{R}S$ is $S-$projective then $M$ is $R-$projective. \\
\end{theorem}

{\bf Proof.} By Lemma \ref{Lemma IV}, $M\otimes_{R}S$ is $S-$projective iff $M\otimes_{R}S'$ is $S'-$projective where $S':=\prod\limits_{\mathfrak{p}\in X}\kappa(\mathfrak{p})$. By the successive applications of Theorems \ref{Corollary VII} and \ref{Corollary VIII}, eventually after finite times we get that $M\otimes_{R}S'$ is $S'-$projective iff $M\otimes_{R}T$ is $T-$projective where $T=\prod\limits_{\mathfrak{p}\in X^{(n)}}\kappa(\mathfrak{p})$. But the kernel of the canonical ring map $R\rightarrow T$ is equal to $\bigcap\limits_{\mathfrak{p}\in X^{(n)}}\mathfrak{p}$. Thus by Theorem \ref{coro 5}, $M$ is $R-$projective. $\Box$ \\

\begin{theorem} Let $\phi:R\rightarrow S$ be a morphism of rings and $X=\Ima\phi^{\ast}$. Suppose there exists some $n\geq1$ such that $\bigcap\limits_{\mathfrak{p}\in X_{(n)}}\mathfrak{p}$ is contained in the radical Jacobson of $R$. If $M\otimes_{R}S$ is $S-$projective then $M$ is $R-$projective. \\
\end{theorem}

{\bf Proof.} It is proven exactly like Theorem \ref{Theorem I}. $\Box$ \\

\section{Auxiliary results}

\begin{proposition} Let $M$ be a finitely generated flat $R-$module, $n$ a natural number and $X=\{\mathfrak{p}\in\Spec(R): \rn_{R_{\mathfrak{p}}}(M_{\mathfrak{p}})=n\}$. Then $\mathcal{Z}(X)=X=\mathcal{F}(X)$. \\
\end{proposition}

{\bf Proof.} By Remark \ref{Remark 040}, we have $X=\Supp N\cap\big(\Spec(R)\setminus\Supp N'\big)$ where $N=\Lambda^{n}(M)$ and $N'=\Lambda^{n+1}(M)$. But $N$ and $N'$ are finitely generated flat $R-$modules and so by Corollary \ref{Corollary I}, the subsets $\Supp(N)=V(\Ann_{R}(N))$ and $\Supp(N')=V(\Ann_{R}(N'))$
are stable under the generalizations. Then the assertion is easily concluded. $\Box$ \\

\begin{proposition}\label{Lemma II} Let $\phi:R\rightarrow S$ be a morphism of rings and $C$ a closed subset of $\Spec(S)$. Then $\overline{\phi^{\ast}(C)}=\mathcal{Z}\big(\phi^{\ast}(C)\big)$. \\
\end{proposition}

{\bf Proof.} Clearly $\mathcal{Z}\big(\phi^{\ast}(C)\big)\subseteq\overline{\phi^{\ast}(C)}$. To see the reverse inclusion, it suffices to show that $\mathcal{Z}\big(\phi^{\ast}(C)\big)$ is closed. There exists an ideal $J$ of $S$ such that $C=V(J)$. We show that $\mathcal{Z}\big(\phi^{\ast}(C)\big)=V(I)$ where $I=\phi^{-1}(J)$. Pick $\mathfrak{p}\in V(I)$. There exists a prime ideal $\mathfrak{p}'$ of $R$ such that it is minimal over $I$ and contained in $\mathfrak{p}$. Thus there exists some $\mathfrak{q}'\in V(J)$ such that $\mathfrak{p}'=\phi^{-1}(\mathfrak{q}')$ since the induced morphism $R/I\rightarrow S/J$ is injective. It follows that $\mathfrak{p}'\in\phi^{\ast}(C)$ and so $\mathfrak{p}\in\mathcal{Z}\big(\phi^{\ast}(C)\big)$. The other direction is obvious. $\Box$ \\

The following result is the dual of Proposition \ref{Lemma II}. \\

\begin{proposition} Let $\phi:R\rightarrow S$ be a morphism of rings and $C$ a flat closed subset of $\Spec(S)$. Then $\overline{\phi^{\ast}(C)}=\mathcal{F}\big(\phi^{\ast}(C)\big)$ where the closure is taken with respect to the flat topology. \\
\end{proposition}

{\bf Proof.} It is proven exactly like Lemma \ref{Lemma II}, recall that the collection of subsets $V(f)$ with $f\in R$ is a sub-basis for the opens of the flat topology, see \cite[Theorem 3.2]{Abolfazl}. $\Box$ \\

\begin{proposition} Let $\phi:R\rightarrow S$ be a morphism of rings and $M$ a finitely generated flat $R-$module with the rank map $\psi$. Then $M\otimes_{R}S$ is a finitely generated flat $S-$module with the rank map $\psi\circ\phi^{\ast}$. \\
\end{proposition}

{\bf Proof.} Let $\theta$ be the rank map of the finitely generated flat  $S-$module $M\otimes_{R}S$ and $\mathfrak{q}$ a prime ideal of $S$ with $\mathfrak{p}=\phi^{-1}(\mathfrak{q})$. Assume that $\theta(\mathfrak{q})=m$ and $\psi(\mathfrak{p})=n$. We have then the following isomorphisms of $S_{\mathfrak{q}}-$modules:
$$(S_{\mathfrak{q}})^{m}\simeq(M\otimes_{R}S)_{\mathfrak{q}}\simeq M_{\mathfrak{p}}\otimes_{R_{\mathfrak{p}}}S_{\mathfrak{q}}\simeq
(S_{\mathfrak{q}})^{n}.$$ This yields that $m=n$, and so $\theta=\psi\circ\phi^{\ast}$. $\Box$ \\

\textbf{Acknowledgements.} The author would likes to give sincere thanks to the referee for very careful reading of the paper and for his/her very valuable and excellent suggestions and comments which greatly improved the paper. We would also like to thank Professors Wolmer Vasconcelos and Jesus Castillo for various scientific correspondences during the writing this paper. \\

\end{document}